\documentclass[12pt]{article}

\usepackage{latexsym,amsfonts,amssymb,mathrsfs,float}
\usepackage{dsfont}
\usepackage{amsthm}
\usepackage[numbers]{natbib}
\usepackage{amsmath}
\usepackage{graphicx}
\usepackage{diagbox}
\usepackage{indentfirst}

\usepackage[dvips]{color}
\usepackage{colordvi,multicol}
\allowdisplaybreaks

\newtheorem{thm}{Theorem}[section]

\newtheorem{pro}[thm]{Proposition}

\numberwithin{equation}{section}

\begin{document}

\title{\bf  A note on the binomial distribution motivated by Chv\'{a}tal's theorem and Tomasewski's theorem}
\author{ Zheng-Yan Guo, Ze-Chun Hu\thanks{Corresponding
author: College of Mathematics, Sichuan University, Chengdu
610065,  China\vskip 0cm E-mail address: zchu@scu.edu.cn}, Run-Yu Wang\\
 {\small College of Mathematics, Sichuan University}}

\date{}

\maketitle

\noindent{\bf Abstract:}\ Let $B(n,p)$ denote a binomial random variable with parameters $n$ and $p$. Chv\'{a}tal's theorem says that for any fixed $n\geq 2$, as $m$ ranges over $\{0,1,\ldots,n\}$, the probability $q_m:=P(B(n,m/n)\leq m)$ is the smallest when  $m$ is closest to $2n/3$. Let $\mathcal{R}$ be the family of random variables of the form $X=\sum^n_{k=1}a_k\varepsilon_k$, where $n\ge 1$,  $a_k, k=1, \dots, n,$ are real numbers with $\sum^n_{k=1} a_k^2=1$, and $\varepsilon_k$, $k=1, 2, \dots$, are independent
Rademacher random variables (i.e., $P(\varepsilon_k=1)=P(\varepsilon_k=-1)=1/2$).
Tomaszewski's theorem says that $\inf_{X\in \mathcal{R}}P(|X|\leq 1)=1/2$. Motivated by
Chv\'{a}tal's Theorem and Tomasewski's Theorem, in this note, we study the minimum value of the probability $f_n(k):=P(|B(n,k/n)-k|\leq \sqrt{{\rm Var} (B(n,k/n))})$ when $k$ ranges over $\{0,1,\ldots,n\}$ for any fixed $n\geq 1$, where ${\rm Var} (\cdot)$ denotes the variance,  and prove that it is the smallest when $k=1$ and $n-1$.

\smallskip


\section{Introduction and main result}

Let $B(n,p)$ denote a binomial random variable with parameters $n$ and $p$. {\it  Chv\'{a}tal's theorem} says that: For any fixed $n\geq 2$, as $m$ ranges over $\{0,1,\ldots,n\}$, the probability $q_m:=P(B(n,m/n)\leq m)$ attains its minimum value when $m$ is closest to $\frac{2n}{3}$. Chv\'{a}tal's theorem may have useful applications, as the probability of a binomial random variable exceeding its mean plays a role in machine learning (see Doerr \cite{Do18}, Greenberg and Mohri \cite{GM14}, Vapnik \cite{Va98} and references therein). For the history of Chv\'{a}tal's theorem, refer to Bababesi et al. \cite{BPR21}, Janson \cite{Ja21} and Sun \cite{Su21}.

 Motivated by Chv\'{a}tal's theorem, a series of papers have studied the infimum value of the probability $P(X\leq E[X])$, where $X$ has the same type of distribution which includes the Poisson, geometric and Pascal distributions (\cite{LXH23}), the negative binomial distribution (\cite{GTH23}), the Gamma distribution (\cite{SHS-a}), the Weibull and Pareto distributions (\cite{LHZ23}),  and some infinitely divisible distributions (\cite{HLZZ24}).

Let $\mathcal{R}$ denote the family of random variables of the form $X=\sum^n_{k=1}a_k\varepsilon_k$, where $n\ge 1$, $a_k, k=1, \dots, n,$ are real numbers satisfying $\sum^n_{k=1} a_k^2=1$, and $\varepsilon_k$, $k=1, 2, \dots$, are independent Rademacher random variables (i.e. $P(\varepsilon_k=1)=P(\varepsilon_k=-1)=1/2$). {\it Tomaszewski's theorem} says that $\inf_{X\in \mathcal{R}}P(|X|\leq 1)=1/2$. Tomaszewski's theorem has wide applications in probability theory, geometric analysis, and computer science. For its history and related problems, see Dvorak and Klein \cite{DK22}, and  Keller and Klein \cite{KK22}.

Motivated by Tomaszewski's theorem, several papers have investigated the infimum value of the probability $P(|X-E[X]|\leq \sqrt{{\rm Var}(X)})$, where the distribution of $X$ belongs to some familiar distributions including the Gamma distribution (\cite{SHS-a}), some infinitely divisible distributions (\cite{SHS-b}),  $F$-distribution (\cite{SHS-c}), and the geometric, symmetric geometric, Poisson and symmetric Poisson distributions (\cite{ZHS25}). If the infimum value is positive for some certain type of distribution, we can say that this kind of distribution possesses the concentration property.

 Motivated by Chv\'{a}tal's theorem and Tomasewski's theorem, in this note, we study the minimum value of the probability $f_n(k):=P(|B(n,k/n)-k|\leq \sqrt{{\rm Var} (B(n,k))})$ for $n\geq 1$ and $k$ ranges over $\{0,1,\ldots,n\}$. The main result is as follows:

\begin{thm}\label{main-thm}
For any fixed positive integer $n$, it holds that
\begin{eqnarray*}
\min_{k=0,1,...,n} f_n(k)=f_n(1)=f_n(n-1)=\left\{
\begin{array}{cl}
1,& \mbox{if}\quad n=1,\\
\left(\frac{n-1}{n}\right)^{n-1},& \mbox{if}\quad n\geq 1.
\end{array}
\right.
\end{eqnarray*}
\end{thm}

\section{Proof}\setcounter{equation}{0}

In this section, we give the proof of Theorem \ref{main-thm}. In Subsection 2.1,  we show that the function $f_n(\cdot)$ possesses a symmetry.  Then in Subsections 2.2, 2.3 and 2.4, we give the proof for  three cases: $n\geq 40,k\geq 10$; $n\geq 40,k<10$; and  $n<40$, respectively. For simplicity, let $X_{n,k}$ denote the random variable $B(n,k)$.

\subsection{Symmetry of the function $f_n(\cdot)$}

\begin{pro}\label{pro-2.1}

For any $n\geq 1$, we have $f_n(k)=f_n(n-k),\forall k=0,1,\ldots,n$.
\end{pro}

\noindent {\bf Proof:} For any $k=0,1,\ldots,n$, we have
\begin{align*}
&{\rm Var}(X_{n,k})=n\cdot \frac{k}{n}\cdot \left(1-\frac{k}{n}\right)=\frac{k(n-k)}{n},\\
&{\rm Var}(X_{n,n-k})=n\cdot \frac{n-k}{n}\cdot \left(1-\frac{n-k}{n}\right)=\frac{k(n-k)}{n}.
\end{align*}
Thus
\begin{align}
&f_n(k)=P\left(|X_{n,k}-k|\leq \sqrt{\frac{k(n-k)}{n}}\right),\label{pro-2.1-a}\\
&f_n(n-k)=P\left(|X_{n,n-k}-(n-k)|\leq \sqrt{\frac{k(n-k)}{n}}\right).\label{pro-2.1-b}
\end{align}
By (\ref{pro-2.1-a}), (\ref{pro-2.1-b}) and the definition of the binomial distribution, we get $f_n(k)=f_n(n-k)$. \hfill\fbox

\subsection{The case of $n\geq 40,k\geq 10$}

By Proposition \ref{pro-2.1}, it is enough to consider the case $k \leq n/2$. Note that
\begin{align*}
f_n(k)&=P\left(|X_{n,k}-E[X_{n,k}]|\leq\sqrt{{\rm Var}(X_{n,k})}\right)=P\left(-1 \leq \frac{X_{n,k}-E[X_{n,k}]}{\sqrt{{\rm Var}(X_{n,k})}} \leq 1\right) .
\end{align*}
Let $Y_{n,k},k=1,\ldots,n$, be independent and identically distributed (i.i.d.) random variables with the distribution  $B(1,k/n)$.  Then $X_{n,k}$ and $\sum_{i=1}^{n} Y_{n,k}$ have the same distribution.  By the  central limit theorem,  it holds that
$$
\frac{X_{n,k}-E[X_{n,k}]}{\sqrt{{\rm Var}(X_{n,k})}} \xrightarrow[]{d} N(0,1)\ \mbox{as}\ n\to \infty.
$$

To complete the proof of this case, we need the convergence rate of the above convergence. Since $Y_{n,k} \sim B(1,k/n)$, we have
\begin{align*}
&E[Y_{n,k}]=\frac{k}{n},\\
&{\rm Var}(Y_{n,k})=1\cdot \frac{k}{n}\cdot \left(1-\frac{k}{n}\right)=\frac{k(n-k)}{n^2},\\
&P\left(|Y_{n,k}-E[Y_{n,k}]|=\frac{k}{n}\right)=1-\frac{k}{n},\\ &P\left(|Y_{n,k}-E[Y_{n,k}]|=1-\frac{k}{n}\right)=\frac{k}{n}.
\end{align*}
It follows that
$$
E[|Y_{n,k}-EY_{n,k}|^3]= \left(\frac{k}{n}\right)^3\left(1-\frac{k}{n}\right)+  \left(1-\frac{k}{n} \right)^3  \frac{k}{n} =\frac{k(n-k)(n^2+2k^2-2nk)}{n^4}.
$$
And thus
\begin{align}
&\sigma:=\sqrt{{\rm Var}(Y_{n,k})}=\frac{\sqrt{k(n-k)}}{n},\nonumber\\
&\rho:=\frac{E[|Y_{n,k}-E[Y_{n,k}]|^3]}{\sigma^3}=\frac{n^2+2k^2-2nk}{n\sqrt{k(n-k)}}.\label{2.3}
\end{align}

By the Berry-Esseen inequality and (\ref{2.3}), we have
\begin{align}\label{2.4}
\sup_{x\in \mathbb{R}}\left|P\left(\frac{X_{n,k}-E[X_{n,k}]}{\sqrt{{\rm Var}(X_{n,k})}} < x\right) - \Phi(x)\right| \leq C_0\frac{\rho}{\sqrt{n}}=
C_0\frac{n^2+2k^2-2nk}{n\sqrt{nk(n-k)}},
\end{align}
where $C_0$ is a positive constant independent of $n$.

By $k \leq n/2$, we know that $n^2+2k^2-2nk=n^2-nk+k(2k-n) \leq n^2-nk$. Then by (\ref{2.4}), we have
\begin{align}\label{2.5}
\sup_{x\in \mathbb{R}}\left|P\left(\frac{X_{n,k}-E[X_{n,k}]}{\sqrt{{\rm Var}(X_{n,k})}} < x\right) - \Phi(x)\right|\leq C_0\frac{n^2-nk}{n\sqrt{nk(n-k)}}=C_0\sqrt{\frac{n-k}{nk}} \leq \frac{C_0}{\sqrt{k}}.
\end{align}

By \cite[Corollary 1]{She11}, we can take $C_0= 0.4748$. Then in virtue of the condition that  $k \geq 10$, we get
\begin{eqnarray*}
\frac{C_0}{\sqrt{k}} \leq \frac{0.4748}{\sqrt{10}}< 0.15014495,
\end{eqnarray*}
which together with (\ref{2.5}) implies that for any $\varepsilon>0$,
\begin{align}
&\left|P\left(\frac{X_{n,k}-E[X_{n,k}]}{\sqrt{{\rm Var}(X_{n,k})}} < -1\right) - \Phi(-1)\right| < 0.15014495,\label{2.6}\\
&\left|P\left(\frac{X_{n,k}-E[X_{n,k}]}{\sqrt{{\rm Var}(X_{n,k})}} < 1+\varepsilon\right) - \Phi(1+\varepsilon)\right|< 0.15014495\label{2.7}
\end{align}
By (\ref{2.6}) and (\ref{2.7}), we get
\begin{eqnarray*}
\left|P\left(-1 \leq \frac{X_{n,k}-E[X_{n,k}]}{\sqrt{{\rm Var}(X_{n,k})}} < 1+\varepsilon\right) - \left(\Phi(1+\varepsilon)-\Phi(-1)\right)\right| < 0.3002899.
\end{eqnarray*}
 Letting $\varepsilon\downarrow 0$, we get
\begin{align}\label{2.8}
\left|P\left(-1 \leq \frac{X_{n,k}-E[X_{n,k}]}{\sqrt{{\rm Var}(X_{n,k})}} \leq 1\right) - \left(\Phi(1)-\Phi(-1)\right)\right|\leq  0.3002899.
\end{align}
Since $\Phi(1)-\Phi(-1) \approx 0.68268949>0.68268948$, by (\ref{2.8}), we know that when $n\geq 40, k \geq 10$,
\begin{eqnarray*}
f_n(k) > 0.68268948 - 0.3002899 =0.38239958.
\end{eqnarray*}

Since $f_n(1)=\left(\frac{n-1}{n}\right)^{n-1}$ is a decreasing function in $n$, and $f_{40}(1)\approx 0.36323244$, we get that for any $n \geq 40, k\geq 10$,
$$
f_n(1)\leq f_{40}(1)< 0.38239958 < f_n(k).
$$

\subsection{The case of $n\geq 40,k< 10$}

The case $k=0$ is trivial. We only need to  prove $f_n(k) \geq f_n(1)$ for  $k=2,3,\ldots,9$. Since $n\geq 40$, it is easy to know that
\begin{align*}
&1<\sqrt{\frac{k(n-k)}{n}}<2,\ \mbox{for}\ k=2,3,4;\\
&2<\sqrt{\frac{k(n-k)}{n}}<3,\ \mbox{for}\ k=5,6,7,8,9.
\end{align*}
It follows that for $k=2,3,4$,
\begin{align*}
\left\{|X_{n,k}-k|\leq \sqrt{\frac{k(n-k)}{n}}\right\}=\{X_{n,k}=k-1,k,k+1\},
\end{align*}
and for $k=5,6,7,8,9$,
\begin{align*}
\left\{|X_{n,k}-k|\leq \sqrt{\frac{k(n-k)}{n}}\right\}=\{X_{n,k}=k-2,k-1,k,k+1,k+2\}.
\end{align*}

In the following, we decompose the proof for $k=2,3,\ldots,9,$ into five cases.

{\it Case 1}: $k=2$. Now
\begin{align*}
f_n(2)&=P(X_{n,2}=1,2,3)\\
&=\frac{1}{n^{n-1}}\left[2(n-2)^{n-1}+\frac{10}{3}(n-1)(n-2)^{n-2}\right].
\end{align*}
Then we have
\begin{align}
&f_n(2)\geq f_n(1)=\left(\frac{n-1}{n}\right)^{n-1}\label{2.9}\\
\Leftrightarrow\ & \frac{1}{n^{n-1}}\left[2(n-2)^{n-1}+\frac{10}{3}(n-1)(n-2)^{n-2}\right]\geq \left(\frac{n-1}{n}\right)^{n-1}\notag\\
\Leftrightarrow\ & 2\left(\frac{n-2}{n-1}\right)^{n-1}+\frac{10}{3}\left(\frac{n-2}{n-1}\right)^{n-2} \geq 1\notag\\
\Leftarrow\ &\frac{10}{3}\left(\frac{n-2}{n-1}\right)^{n-2} \geq 1.\notag
\end{align}
Since  $\left(\frac{x-2}{x-1}\right)^{x-2}$ is  a decreasing function of $x$ when $x \geq 40$, and  the limit is $\frac{1}{e}$, we have that  $\frac{10}{3}\left(\frac{n-2}{n-1}\right)^{n-2}
\geq \frac{10}{3e}> 1$. Hence (\ref{2.9}) holds.
\smallskip

{\it Case 2}: $k=3$. Now
\begin{align*}
f_n(3)&=P(X_{n,3}=2,3,4)\\
&=\frac{1}{n^{n-1}}\left[\frac{9}{2}(n-1)(n-3)^{n-2}+\frac{63}{8}(n-1)(n-2)(n-3)^{n-3}\right].
\end{align*}
Then we have
\begin{align}
& f_n(3) \geq f_n(1)=\left(\frac{n-1}{n}\right)^{n-1}\label{2.10}\\
\Leftrightarrow\ &\frac{1}{n^{n-1}}\left[\frac{9}{2}(n-1)(n-3)^{n-2}+\frac{63}{8}(n-1)(n-2)(n-3)^{n-3}\right]
\geq \left(\frac{n-1}{n}\right)^{n-1}\notag\\
\Leftrightarrow\ & \frac{9}{2}\left(1-\frac{2}{n-1}\right)^{n-2}+\frac{63}{8} \left(1-\frac{1}{n-1}\right)\left(1-\frac{2}{n-1}\right)^{n-3} \geq 1.\notag
\end{align}
It is easy to know that $\left(1-\frac{2}{x-1}\right)^{x-2}$ and $1-
\frac{1}{x-1}$ are two  increasing functions, $\left(1-\frac{2}{x-1}\right)^{x-3}$ is a decreasing function of $x$ on $[40,\infty)$, and  $\lim\limits_{x\to\infty}\left(1-\frac{2}{x-1}\right)^{x-3}=\frac{1}{e^2}$.
It follows that
\begin{align*}
&\frac{9}{2}\left(1-\frac{2}{n-1}\right)^{n-2}+\frac{63}{8} \left(1-\frac{1}{n-1}\right)\left(1-\frac{2}{n-1}\right)^{n-3}\\
&\geq \frac{9}{2} \left(\frac{37}{39}\right)^{38}+\frac{63}{8}\cdot \frac{38}{39} \cdot \frac{1}{e^2}\approx 1.65>1.
\end{align*}
Hence (\ref{2.10}) holds.
\smallskip

{\it Case 3}: $k=4$. Now
\begin{align*}
f_n(4)&=P(X_{n,4}=3,4,5)\\
&=\frac{1}{n^{n-1}} \left[\frac{32}{3}(n-1)(n-2)(n-4)^{n-3}+\frac{96}{5}(n-1)(n-2)(n-3)(n-4)^{n-4}\right].
\end{align*}
Then we have
\begin{align}
&f_n(4)\geq f_n(1)=\left(\frac{n-1}{n}\right)^{n-1}\label{2.11}\\
\Leftrightarrow \ &\frac{1}{n^{n-1}} \left[\frac{32}{3}(n-1)(n-2)(n-4)^{n-3}+\frac{96}{5}(n-1)(n-2)(n-3)(n-4)^{n-4}\right]\geq \left(\frac{n-1}{n}\right)^{n-1}\notag \\
\Leftrightarrow\ &\frac{32}{3}\left(1-\frac{1}{n-1}\right)\left(1-\frac{3}{n-1}\right)^{n-3}+
\frac{96}{5}\left(1-\frac{1}{n-1}\right)\left(1-\frac{2}{n-1}\right)
\left(1-\frac{3}{n-1}\right)^{n-4}\geq 1\notag\\
\Leftarrow\ &\left(\frac{32}{3}+\frac{96}{5}\right)\left(1-\frac{3}{n-1}\right)^{n-2}\geq 1.\notag
\end{align}
Since $\left(1-\frac{3}{x-1}\right)^{x-2}$ is an increasing function of $x$ on $[40,\infty)$, we have
\begin{align*}
\left(\frac{32}{3}+\frac{96}{5}\right)\left(1-\frac{3}{n-1}\right)^{n-2}
\geq \left(\frac{32}{3}+\frac{96}{5}\right)\left(\frac{36}{39}\right)^{38} \approx 1.42635>1.
\end{align*}
Hence (\ref{2.11}) holds.

\smallskip

{\it Case 4}: $k=5,6,7,8$. Now
\begin{align*}
f_n(k)&=P(X_{n,k}=k-2,k-1,k,k+1,k+2)\\
&=\sum_{i=k-2}^{k+2}\binom{n}{i}\left(\frac{k}{n}\right)^i\left(1-\frac{k}{n}\right)^{n-i}\\
&=\sum_{i=k-2}^{k+2}\frac{n!}{i!(n-i)!}\left(\frac{k}{n}\right)^i
\left(1-\frac{k}{n}\right)^{n-i}.
\end{align*}
Then we have
\begin{align}
&f_n(k)\geq f_n(1)=\left(\frac{n-1}{n}\right)^{n-1}\label{2.12}\\
\Leftrightarrow\ &\sum_{i=k-2}^{k+2}\frac{n!}{i!(n-i)!}\left(\frac{k}{n}\right)^i
\left(1-\frac{k}{n}\right)^{n-i}\geq \left(\frac{n-1}{n}\right)^{n-1}\notag\\
\Leftrightarrow\ &\sum_{i=k-2}^{k+2}\frac{k^i}{i!}\frac{(n-2)\cdots(n-i+1)(n-k)^{n-i}}{(n-1)^{n-2}}\geq 1\notag\\
\Leftarrow\ &\sum_{i=k-2}^{k+2}\frac{k^i}{i!}\left(1-\frac{k-1}{n-1}\right)^{n-2}\geq 1\notag\\
\Leftarrow\ &\sum_{i=k-2}^{k+2}\frac{k^i}{i!}\left(1-\frac{k-1}{40-1}\right)^{40-2}\geq 1,
\notag
\end{align}
where in the last step, we used the fact that $\left(1-\frac{k-1}{x-1}\right)^{x-2}$ is an increasing function of $x$ on $[40,\infty)$ for any $k=5,6,7,8$.

Define
$$
C_k:=\sum_{i=k-2}^{k+2}\frac{k^i}{i!}\left(1-\frac{k-1}{40-1}\right)^{40-2}.
$$
By calculation, we get
$$
C_5\approx 1.80299, \ C_6\approx 1.52806,\ C_7\approx 1.26193,\ C_8\approx 1.01213.
$$
Hence (\ref{2.12}) holds.
\smallskip

{\it Case 5}: $k=9$. Now we have
\begin{align}\label{2.13}
f_n(9)&=P(X_{n,9}=7,8,9,10,11)\nonumber\\
&=\sum_{i=7}^{11}\frac{n!}{i!(n-i)!}\left(\frac{9}{n}\right)^i
\left(1-\frac{9}{n}\right)^{n-i}.
\end{align}
Then we have
\begin{align}
&f_n(9)\geq f_n(1)=\left(\frac{n-1}{n}\right)^{n-1}\label{2.14}\\
\Leftrightarrow\ &\sum_{i=7}^{9}\frac{n!}{i!(n-i)!}\left(\frac{9}{n}\right)^i
\left(1-\frac{9}{n}\right)^{n-i}\geq \left(\frac{n-1}{n}\right)^{n-1}\notag\\
\Leftrightarrow\ &\sum_{i=7}^{11}\frac{9^i}{i!}\frac{(n-2)(n-3)\cdots(n-i+1)(n-9)^{n-i}}{(n-1)^{n-2}}\geq 1\notag\\
\Leftrightarrow\ &\sum_{i=7}^{11}\frac{9^i}{i!}\left(1-\frac{8}{n-1}\right)^{n-2}
\cdot \frac{(n-2)(n-3)\cdots(n-i+1)}{(n-9)^{i-2}}\geq 1\notag\\
\Leftarrow\ &\sum_{i=7}^{9}\frac{9^i}{i!}\left(1-\frac{8}{n-1}\right)^{n-2}\geq 1,\label{2.15}
\end{align}
where in the last step, we used the fact that
\begin{align}\label{2.16}
\frac{(n-2)(n-3)\cdots(n-i+1)}{(n-9)^{i-2}}\geq 1,\quad \forall i=7,\cdots,11.
\end{align}
For $i=7,8,9,10$, we have  $n-i+1\geq n-9$, and thus (\ref{2.16}) holds. For $i=11$, we have
\begin{align*}
&\frac{(n-2)(n-3)\cdots(n-10)}{(n-9)^9}\geq 1\\
\Leftarrow\ &\frac{(n-2)(n-10)}{(n-9)^2}\geq 1\\
\Leftrightarrow\ &6n\geq 61\\
\Leftarrow\ & n\geq 40.
\end{align*}
Hence (\ref{2.16}) is true.

Since $(1-\frac{8}{x-1})^{x-2}$ is an increasing function of $x$ on $[100,\infty)$, we get that for any $n\geq 100$,
$$
\sum_{i=7}^{9}\frac{9^i}{i!}\left(1-\frac{8}{n-1}\right)^{n-2}\geq \sum_{i=7}^{9}\frac{9^i}{i!}\left(1-\frac{8}{99}\right)^{98}\approx 1.25277>1.
$$
It follows that (\ref{2.15}) holds and thus (\ref{2.14}) holds for any $n\geq 100$.

For any $n=40,41,\cdots,99$, by direct calculations in virtue of (\ref{2.13}) and $f_n(1)=\left(\frac{n-1}{n}\right)^{n-1}$, we get that $$
f_n(9)>0.61>0.38>f_n(1).
$$
In a word, (\ref{2.14}) holds for any $n\geq 40$.

\subsection{The case of $n<40$}

In this case, we directly calculate the numerical results of $f_n(k)$ for $n=1,2,\cdots,39$ and $f=1,2,\cdots,\lfloor\frac{n}{2}\rfloor$ by virtue of Proposition \ref{pro-2.1}. The case that $n=1,2$ is trivial. For $n=39, k=1,2,\cdots,19$, the numerical results of $f_{39}(k)$ are as follows:
\begin{align*}
&0.372668, 0.733642, 0.634857, 0.571572, 0.773351, 0.735495, 0.704738, 0.679446, 0.658464, 0.640964, \\
&0.62634, 0.614139, 0.604025, 0.595739, 0.751145, 0.745834, 0.741918, 0.739339, 0.738058.
\end{align*}
It follows that $f_{39}(1)$ is the smallest. For the case that $n=3,4,\cdots,38$, the numerical results are presented in Table 1 and Table 2, which together with Proposition \ref{pro-2.1} and the fact that $f_n(0)=f_n(n)=1$ show that $f_n(1)=f_n(n-1)$ is the smallest in the set  $\{f_n(k): k=0,1,2,\cdots,n\}$.

\begin{table}
\centering
\tiny
\begin{tabular}{|l|l|l|l|l|l|l|l|l|l|l|l|l|l|l|l|l|l|l|}
\hline
\diagbox{$k$}{$n$} & 3     & 4     & 5     & 6     & 7     & 8     & 9     & 10    & 11    & 12    & 13    & 14    & 15    & 16    & 17    & 18    & 19    & 20     \\
\hline
1                & 0.44~ & 0.42~ & 0.41~ & 0.40~ & 0.40~ & 0.39~ & 0.39~ & 0.39~ & 0.39~ & 0.38~ & 0.38~ & 0.38~ & 0.38~ & 0.38~ & 0.38~ & 0.38~ & 0.38~ & 0.38~  \\
\hline
2                & 0.44~ & 0.88~ & 0.84~ & 0.81~ & 0.80~ & 0.79~ & 0.78~ & 0.77~ & 0.77~ & 0.76~ & 0.76~ & 0.76~ & 0.75~ & 0.75~ & 0.75~ & 0.75~ & 0.75~ & 0.75~  \\
\hline
3                &       & 0.42~ & 0.84~ & 0.78~ & 0.75~ & 0.73~ & 0.71~ & 0.70~ & 0.69~ & 0.68~ & 0.68~ & 0.67~ & 0.67~ & 0.66~ & 0.66~ & 0.66~ & 0.66~ & 0.65~  \\
\hline
4                &       &       & 0.41~ & 0.81~ & 0.75~ & 0.71~ & 0.69~ & 0.67~ & 0.65~ & 0.64~ & 0.63~ & 0.62~ & 0.62~ & 0.61~ & 0.61~ & 0.60~ & 0.60~ & 0.60~  \\
\hline
5                &       &       &       & 0.40~ & 0.80~ & 0.73~ & 0.69~ & 0.66~ & 0.64~ & 0.62~ & 0.61~ & 0.60~ & 0.59~ & 0.58~ & 0.57~ & 0.57~ & 0.56~ & 0.56~  \\
\hline
6                &       &       &       &       & 0.40~ & 0.79~ & 0.71~ & 0.67~ & 0.64~ & 0.61~ & 0.59~ & 0.58~ & 0.57~ & 0.56~ & 0.55~ & 0.79~ & 0.78~ & 0.78~  \\
\hline
7                &       &       &       &       &       & 0.39~ & 0.78~ & 0.70~ & 0.65~ & 0.62~ & 0.59~ & 0.58~ & 0.56~ & 0.55~ & 0.78~ & 0.77~ & 0.77~ & 0.76~  \\
\hline
8                &       &       &       &       &       &       & 0.39~ & 0.77~ & 0.69~ & 0.64~ & 0.61~ & 0.58~ & 0.56~ & 0.79~ & 0.78~ & 0.77~ & 0.76~ & 0.75~  \\
\hline
9                &       &       &       &       &       &       &       & 0.39~ & 0.77~ & 0.68~ & 0.63~ & 0.60~ & 0.57~ & 0.55~ & 0.78~ & 0.76~ & 0.75~ & 0.74~  \\
\hline
10               &       &       &       &       &       &       &       &       & 0.39~ & 0.76~ & 0.68~ & 0.62~ & 0.59~ & 0.56~ & 0.78~ & 0.77~ & 0.75~ & 0.74~  \\
\hline
11               &       &       &       &       &       &       &       &       &       & 0.38~ & 0.76~ & 0.67~ & 0.62~ & 0.58~ & 0.55~ & 0.77~ & 0.76~ & 0.74~  \\
\hline
12               &       &       &       &       &       &       &       &       &       &       & 0.38~ & 0.76~ & 0.67~ & 0.61~ & 0.57~ & 0.79~ & 0.77~ & 0.75~  \\
\hline
13               &       &       &       &       &       &       &       &       &       &       &       & 0.38~ & 0.75~ & 0.66~ & 0.61~ & 0.57~ & 0.78~ & 0.76~  \\
\hline
14               &       &       &       &       &       &       &       &       &       &       &       &       & 0.38~ & 0.75~ & 0.66~ & 0.60~ & 0.56~ & 0.78~  \\
\hline
15               &       &       &       &       &       &       &       &       &       &       &       &       &       & 0.38~ & 0.75~ & 0.66~ & 0.60~ & 0.56~  \\
\hline
16               &       &       &       &       &       &       &       &       &       &       &       &       &       &       & 0.38~ & 0.75~ & 0.66~ & 0.60~  \\
\hline
17               &       &       &       &       &       &       &       &       &       &       &       &       &       &       &       & 0.38~ & 0.75~ & 0.65~  \\
\hline
18               &       &       &       &       &       &       &       &       &       &       &       &       &       &       &       &       & 0.38~ & 0.75~  \\
\hline
19               &       &       &       &       &       &       &       &       &       &       &       &       &       &       &       &       &       & 0.38~  \\
\hline
20               &       &       &       &       &       &       &       &       &       &       &       &       &       &       &       &       &       &        \\
\hline
\end{tabular}
\smallskip

\centerline{Tabel 1}

\end{table}

\begin{table}
\centering
\tiny
\begin{tabular}{|l|l|l|l|l|l|l|l|l|l|l|l|l|l|l|l|l|l|l|}
\hline
\diagbox{$k$}{$n$} & 21    & 22    & 23    & 24    & 25    & 26    & 27    & 28    & 29    & 30    & 31    & 32    & 33    & 34    & 35    & 36    & 37    & 38     \\
\hline
1                & 0.38~ & 0.38~ & 0.38~ & 0.38~ & 0.38~ & 0.38~ & 0.37~ & 0.37~ & 0.37~ & 0.37~ & 0.37~ & 0.37~ & 0.37~ & 0.37~ & 0.37~ & 0.37~ & 0.37~ & 0.37~  \\
\hline
2                & 0.74~ & 0.74~ & 0.74~ & 0.74~ & 0.74~ & 0.74~ & 0.74~ & 0.74~ & 0.74~ & 0.74~ & 0.74~ & 0.74~ & 0.74~ & 0.74~ & 0.74~ & 0.73~ & 0.73~ & 0.73~  \\
\hline
3                & 0.65~ & 0.65~ & 0.65~ & 0.65~ & 0.65~ & 0.64~ & 0.64~ & 0.64~ & 0.64~ & 0.64~ & 0.64~ & 0.64~ & 0.64~ & 0.64~ & 0.64~ & 0.64~ & 0.64~ & 0.64~  \\
\hline
4                & 0.60~ & 0.59~ & 0.59~ & 0.59~ & 0.59~ & 0.59~ & 0.58~ & 0.58~ & 0.58~ & 0.58~ & 0.58~ & 0.58~ & 0.58~ & 0.58~ & 0.57~ & 0.57~ & 0.57~ & 0.57~  \\
\hline
5                & 0.56~ & 0.55~ & 0.55~ & 0.55~ & 0.79~ & 0.79~ & 0.79~ & 0.79~ & 0.79~ & 0.78~ & 0.78~ & 0.78~ & 0.78~ & 0.78~ & 0.78~ & 0.78~ & 0.78~ & 0.77~  \\
\hline
6                & 0.77~ & 0.77~ & 0.77~ & 0.76~ & 0.76~ & 0.76~ & 0.76~ & 0.75~ & 0.75~ & 0.75~ & 0.75~ & 0.75~ & 0.74~ & 0.74~ & 0.74~ & 0.74~ & 0.74~ & 0.74~  \\
\hline
7                & 0.75~ & 0.75~ & 0.74~ & 0.74~ & 0.74~ & 0.73~ & 0.73~ & 0.73~ & 0.72~ & 0.72~ & 0.72~ & 0.72~ & 0.71~ & 0.71~ & 0.71~ & 0.71~ & 0.71~ & 0.71~  \\
\hline
8                & 0.74~ & 0.73~ & 0.73~ & 0.72~ & 0.72~ & 0.71~ & 0.71~ & 0.71~ & 0.70~ & 0.70~ & 0.70~ & 0.69~ & 0.69~ & 0.69~ & 0.69~ & 0.68~ & 0.68~ & 0.68~  \\
\hline
9                & 0.73~ & 0.72~ & 0.71~ & 0.71~ & 0.70~ & 0.70~ & 0.69~ & 0.69~ & 0.68~ & 0.68~ & 0.68~ & 0.67~ & 0.67~ & 0.67~ & 0.67~ & 0.66~ & 0.66~ & 0.66~  \\
\hline
10               & 0.73~ & 0.72~ & 0.71~ & 0.70~ & 0.69~ & 0.69~ & 0.68~ & 0.68~ & 0.67~ & 0.67~ & 0.66~ & 0.66~ & 0.66~ & 0.65~ & 0.65~ & 0.65~ & 0.65~ & 0.64~  \\
\hline
11               & 0.73~ & 0.71~ & 0.70~ & 0.69~ & 0.69~ & 0.68~ & 0.67~ & 0.67~ & 0.66~ & 0.66~ & 0.65~ & 0.65~ & 0.64~ & 0.64~ & 0.64~ & 0.63~ & 0.63~ & 0.63~  \\
\hline
12               & 0.73~ & 0.72~ & 0.70~ & 0.69~ & 0.68~ & 0.67~ & 0.67~ & 0.66~ & 0.65~ & 0.65~ & 0.64~ & 0.64~ & 0.63~ & 0.63~ & 0.63~ & 0.62~ & 0.62~ & 0.62~  \\
\hline
13               & 0.74~ & 0.72~ & 0.71~ & 0.69~ & 0.68~ & 0.67~ & 0.66~ & 0.66~ & 0.65~ & 0.64~ & 0.64~ & 0.63~ & 0.63~ & 0.62~ & 0.62~ & 0.61~ & 0.61~ & 0.61~  \\
\hline
14               & 0.75~ & 0.73~ & 0.71~ & 0.70~ & 0.69~ & 0.67~ & 0.66~ & 0.66~ & 0.65~ & 0.64~ & 0.63~ & 0.63~ & 0.62~ & 0.62~ & 0.61~ & 0.61~ & 0.60~ & 0.60~  \\
\hline
15               & 0.77~ & 0.75~ & 0.73~ & 0.71~ & 0.69~ & 0.68~ & 0.67~ & 0.66~ & 0.65~ & 0.64~ & 0.63~ & 0.62~ & 0.62~ & 0.61~ & 0.61~ & 0.60~ & 0.60~ & 0.76~  \\
\hline
16               & 0.56~ & 0.77~ & 0.74~ & 0.72~ & 0.70~ & 0.69~ & 0.67~ & 0.66~ & 0.65~ & 0.64~ & 0.63~ & 0.62~ & 0.62~ & 0.61~ & 0.60~ & 0.60~ & 0.75~ & 0.75~  \\
\hline
17               & 0.60~ & 0.55~ & 0.77~ & 0.74~ & 0.72~ & 0.70~ & 0.68~ & 0.67~ & 0.65~ & 0.64~ & 0.63~ & 0.62~ & 0.62~ & 0.61~ & 0.60~ & 0.60~ & 0.75~ & 0.75~  \\
\hline
18               & 0.65~ & 0.59~ & 0.55~ & 0.76~ & 0.74~ & 0.71~ & 0.69~ & 0.68~ & 0.66~ & 0.65~ & 0.64~ & 0.63~ & 0.62~ & 0.61~ & 0.60~ & 0.76~ & 0.75~ & 0.74~  \\
\hline
19               & 0.74~ & 0.65~ & 0.59~ & 0.55~ & 0.76~ & 0.73~ & 0.71~ & 0.69~ & 0.67~ & 0.66~ & 0.64~ & 0.63~ & 0.62~ & 0.61~ & 0.60~ & 0.60~ & 0.75~ & 0.74~  \\
\hline
20               & 0.38~ & 0.74~ & 0.65~ & 0.59~ & 0.79~ & 0.76~ & 0.73~ & 0.71~ & 0.68~ & 0.67~ & 0.65~ & 0.64~ & 0.63~ & 0.62~ & 0.61~ & 0.60~ & 0.75~ & 0.74~  \\
\hline
21               &       & 0.38~ & 0.74~ & 0.65~ & 0.59~ & 0.79~ & 0.76~ & 0.73~ & 0.70~ & 0.68~ & 0.66~ & 0.65~ & 0.63~ & 0.62~ & 0.61~ & 0.60~ & 0.75~ & 0.75~  \\
\hline
22               &       &       & 0.38~ & 0.74~ & 0.65~ & 0.59~ & 0.79~ & 0.75~ & 0.72~ & 0.70~ & 0.68~ & 0.66~ & 0.64~ & 0.63~ & 0.62~ & 0.61~ & 0.60~ & 0.75~  \\
\hline
23               &       &       &       & 0.38~ & 0.74~ & 0.64~ & 0.58~ & 0.79~ & 0.75~ & 0.72~ & 0.70~ & 0.67~ & 0.66~ & 0.64~ & 0.63~ & 0.61~ & 0.60~ & 0.76~  \\
\hline
24               &       &       &       &       & 0.38~ & 0.74~ & 0.64~ & 0.58~ & 0.79~ & 0.75~ & 0.72~ & 0.69~ & 0.67~ & 0.65~ & 0.64~ & 0.62~ & 0.61~ & 0.60~  \\
\hline
25               &       &       &       &       &       & 0.38~ & 0.74~ & 0.64~ & 0.58~ & 0.78~ & 0.75~ & 0.72~ & 0.69~ & 0.67~ & 0.65~ & 0.63~ & 0.62~ & 0.61~  \\
\hline
26               &       &       &       &       &       &       & 0.37~ & 0.74~ & 0.64~ & 0.58~ & 0.78~ & 0.75~ & 0.71~ & 0.69~ & 0.67~ & 0.65~ & 0.63~ & 0.62~  \\
\hline
27               &       &       &       &       &       &       &       & 0.37~ & 0.74~ & 0.64~ & 0.58~ & 0.78~ & 0.74~ & 0.71~ & 0.69~ & 0.66~ & 0.65~ & 0.63~  \\
\hline
28               &       &       &       &       &       &       &       &       & 0.37~ & 0.74~ & 0.64~ & 0.58~ & 0.78~ & 0.74~ & 0.71~ & 0.68~ & 0.66~ & 0.64~  \\
\hline
29               &       &       &       &       &       &       &       &       &       & 0.37~ & 0.74~ & 0.64~ & 0.58~ & 0.78~ & 0.74~ & 0.71~ & 0.68~ & 0.66~  \\
\hline
30               &       &       &       &       &       &       &       &       &       &       & 0.37~ & 0.74~ & 0.64~ & 0.58~ & 0.78~ & 0.74~ & 0.71~ & 0.68~  \\
\hline
31               &       &       &       &       &       &       &       &       &       &       &       & 0.37~ & 0.74~ & 0.64~ & 0.57~ & 0.78~ & 0.74~ & 0.71~  \\
\hline
32               &       &       &       &       &       &       &       &       &       &       &       &       & 0.37~ & 0.74~ & 0.64~ & 0.57~ & 0.78~ & 0.74~  \\
\hline
33               &       &       &       &       &       &       &       &       &       &       &       &       &       & 0.37~ & 0.74~ & 0.64~ & 0.57~ & 0.77~  \\
\hline
34               &       &       &       &       &       &       &       &       &       &       &       &       &       &       & 0.37~ & 0.73~ & 0.64~ & 0.57~  \\
\hline
35               &       &       &       &       &       &       &       &       &       &       &       &       &       &       &       & 0.37~ & 0.73~ & 0.64~  \\
\hline
36               &       &       &       &       &       &       &       &       &       &       &       &       &       &       &       &       & 0.37~ & 0.73~  \\
\hline
37               &       &       &       &       &       &       &       &       &       &       &       &       &       &       &       &       &       & 0.37~  \\
\hline
\end{tabular}

\smallskip
\centerline{Tabel 2}
\end{table}

\newpage

\noindent {\bf\large Acknowledgments}\quad  This work was supported by the National Natural Science Foundation of China (12171335).

\end{document}